\documentclass[12pt]{article}
\usepackage{makeidx,amsfonts,graphicx,amsthm, amsmath}
\newtheorem{theorem}{Theorem}{\bfseries}{\rmfamily}
\newtheorem{lemma}[theorem]{Lemma}{\bfseries}{\rmfamily}
{\bfseries}{\rmfamily}
{\bfseries}{\rmfamily}
\newtheorem{conjecture}[theorem]{Conjecture}{\bfseries}{\rmfamily}
\def\rr{{\mathbb R}}

\def\nn{{\mathbb N}}

\begin{document}
\centerline{\Large\bf On an Improvement of a Result by}

\vspace*{ 3 mm} 

\centerline{\Large\bf Niederreiter and Wang
Concerning the}

\vspace*{ 3 mm} 

\centerline{\Large\bf Expected Linear Complexity of Multisequences}
\vskip+0.5cm

\footnotetext{Supported by grants FONDECYT 7060126, RFFI 06-01-00518,
MD-3003.2006.1 (N.M.) and FONDECYT  1040975  (M.V.)}

\vskip+0.5cm
 \centerline{\bf
Moshchevitin N.${}^\dag$ \hspace*{ 3 cm} Vielhaber M.${}^\ddag$
 }

\vskip+0.5cm

\begin{minipage}[t]{7 cm}
${}^\dag$
Department of Theory of Numbers\\
Fac.~Mathematics and Mechanics\\
Moscow State University\\
Leninskie Gory\\
119992, Moscow\\
RUSSIA\\
{\tt moshchevitin@rambler.ru}
\end{minipage}
\begin{minipage}[t]{2 cm}
\end{minipage}
\begin{minipage}[t]{8 cm}
${}^\ddag$ Instituto de Matem\'aticas\\
Facultad de Ciencias\\
Universidad Austral de Chile\\
Casilla 567\\
Valdivia\\
CHILE\\
{\tt vielhaber@gmail.com}
\end{minipage}
\hfill\hspace*{1 mm}

  \vskip+1.0cm

\begin{quote}
{\bf Abstract.}

We show that the expected value for the joint linear complexity of an $m$--multisequence of length $n$ is
$$E_n^{(m)}=\left\lceil n\frac{m}{m+1}\right\rceil+O(1),$$
improving on a result by Niederreiter and Wang.
\end{quote}

\noindent MSC: 11B85 Automata sequences
\\\\
Let $T= ( \{s_1 \}_{j=1}^{n},...,  \{s_m \}_{j=1}^{n}) \in (F_q^m)^n $ denote a multisequence and let $L_n^{(m)} (T)$ denote its joint linear
complexity. Let $N_n^{(m)}(L)$ denote the number of multisequences $ $ with $L_n^{(m)} (T)= L$ and let
$$
 E_n^{(m)}
=
q^{-nm}
\sum_{T
\in (F_q^{(m)})^n}
L_n^{(m)} (T)
=
q^{-nm}
\sum_{L=0}^n
L
N_n^{(m)}(L).
$$
be the expected joint linear complexity for a prefix of length $n$ of an $m$--multisequence.

Rueppel \cite{RU}
proved that
$$
 E_n^{(1)}
=\frac{n}{2} + O(1).
$$

Niederreiter and Wang
proved
$$
 E_n^{(2)}
=\frac{2n}{3} + O(1)
\ \ \text{and}\ \
 E_n^{(3)}
=\frac{3n}{4} + O(1)
$$
in \cite{NW1} and \cite{NW3}, respectively,
and in \cite{NW2} they obtain
$$
 E_n^{(m)}
=\frac{mn}{m+1} + o(n)
$$
for any natural $m$.

We improve the details of the method from \cite{NW1},
\cite{NW2} and establish the following result.

\begin{theorem}
For all $m\in\nn$,
$$
 E_n^{(m)}
=\frac{mn}{m+1} + O(1).
$$
\end{theorem}

Here the constant in the symbol $O(1)$ may depend on $q$ and $m$ but {\it does not} depend on $n$.

The proof will follow from the powerful formula for the value
$N_n^{(m)}(L)$ in terms of partitions of  $L$ into $M$ parts
given in \cite[formula 2]{NW2} or \cite[formula (11) and Theorem
  2]{NW1}.

We do not need  to describe Niederreiter's and Wang's formula in
the exact way,
we need only two corollaries of the formula.
These corollaries are presented in \cite{NW2}.

\begin{lemma} {\rm \cite[Lemma 2]{NW2}}
$$
N_n^{(m)}(L)
\le q^{(m+1)L}.
$$
\end{lemma}

\begin{lemma} {\rm \cite[Lemma 3]{NW2}}
Let $ i_1\ge i_2 \ge ...\ge i_m\ge 0$
be integers with $i_1+...+i_m = L$ and
$P(m,L)$
denote the set of all m-tuples $I = (i_1, i_2 , ..., i_m)$
under the condition specified.
Then
$$
N_n^{(m)}(L) \le c(q,m) \sum_{I\in P(m,L)} q^{2\sum_{k=1}^m (k-1)i_k+ 2m(n-L)}
$$
where $c(q,m)$ depends only on $q$ and $m$.
\end{lemma}

In fact this lemma follows from the proof of Lemma 3 from
\cite{NW2}, formulas (3--5).

We first deduce a corollary from Lemma 2. However,  instead of using arguments
from \cite[Lemma 1]{NW2}
we use an estimation on integer points in a polytope--type domain.

We remind that Lemma 1 from \cite{NW2} establishes for reals
$x_1\ge x_2 \ge ...\ge x_m\ge 0$ such that  $x_1+...+x_m = L$
the inequality
$$
2
\sum_{k=1}^m (k-1)x_k
\le
(m-1)
\sum_{k=1}^m x_k =(m-1)L
.
$$

It means that if we consider the set
$$
\Omega =
\{
 x = (x_1,...,x_m)\in \rr^m:x_1\ge x_2 \ge ...\ge x_m\ge 0, x_1+...+x_m = L\}
$$
then
$$
\max_{x\in\Omega}2
\sum_{k=1}^m (k-1)x_k =(m-1)L,
$$
   In fact $\Omega$ is an 
    $m-1$-dimensional simplex in
$\mathbb{R}^m$. Let  $1\le \nu \le m$. Denote by $x^\nu \in \mathbb{R}^m$ the point
whose first $\nu$ coordinates are equal to
$\frac{L}{\nu} $ and all other coordinates are equal to zero. Then $x^\nu , 1\le \nu \le m$, 
 are the all vertices of the simplex $\Omega$.
  The linear function $2\sum_{k=1}^m (k-1)x_k$ attains its maximum on
  some vertex, it is easily
verified that this point is unique.

Note that from Lemma 3 it follows that
\begin{equation}
N_n^{(m)}(L)
\le c(q,m)
\sum_{H=0}^{(m-1)L}\rho_H
q^{2mn -(m+1)L -H},
\label{a}
\end{equation}
where $\rho_H$ is the number of integer  solutions of the system
$$
\begin{cases}
 i_1\ge i_2 \ge ...\ge i_m\ge 0,
\\
i_1+...+i_m =L,
\\
2\sum_{k=1}^m (k-1)i_k = (m-1)L-H.
\end{cases}
$$

Obviously $\rho_H = 0$ for $H\not\equiv(m-1)L\mod 2$, but this is not
of importance.

Let
$$
\Omega_H =
\Omega \cap\left\{ 
 x = (x_1,...,x_m)\in R^m:
2\sum_{k=1}^m (k-1)x_k
\ge (m-1)L - H
\right\} 
$$
be an 
 $m-1$-dimensional polytope in $\mathbb{R}^m$ and $M_H$ be the number of integer points in $\Omega_H$.
  Obviously $ \rho_H \le M_H$.
If $ H =0$ then $\Omega_0$ contains  only the point $x^m$. We may see that $ \Omega_H $ consists of all points $x=(x_1,...,x_m)$ such that
$$
\begin{cases}
 x_1\ge x_2 \ge ...\ge x_m\ge 0,
\\
x_1+...+x_m =L,
\\
2\sum_{k=1}^m (k-1)x_k \ge (m-1)L-H.
\end{cases}
$$

\begin{lemma}
The   number $M_H$ of integer points in
$\Omega_H, H \ge 0$
can be bounded by
\begin{equation}
M_H \le  (H+1)^m. \label{b}
\end{equation}

\end{lemma}

\begin{proof}
Instead of the polytope $ \Omega_H $ we consider the polytope $ \Omega_H ^* $ defined by the conditions
$$
\begin{cases}
 x_1\ge x_2 \ge ...\ge x_m,
\\
x_1+...+x_m =L,
\\
2\sum_{k=1}^m (k-1)x_k \ge (m-1)L-H.
\end{cases}
$$
Here the inequality $ x_m \ge 0$ is omitted and hence $ \Omega_H \subseteq \Omega_H^*$.

Now we  see that the polytope  $\Omega_H^*$ is an 
$m-1$-dimensional simplex with vertices
 $x(H,\nu ), \nu = 1, ..., m$.
    These
vertices
 can be easily calculated:
We 
have  $x(H,m) = x^m$, 
 and for $1\le \nu \le m-1$, 
let $t_\nu = \frac{H}{(m-\nu)L}\ge 0$ and then
$$
x(H,\nu ) = x^m(1-t_\nu) + x^{\nu} t_\nu. 
$$

To prove 
this statement we must consider
the intersection of the hyperplane
$$
\left\{(x_1,...,x_m) \in\mathbb{R}^m:\,\,\, 2\sum_{k=1}^m (k-1)x_k = (m-1)L-H\right\}
$$
with the straight lines containing the edges $ [x^m,x^{\nu }], \nu = 1,...,m-1$ of the simplex $\Omega$.  The 
hyperplane intersects  all these
lines, namely the rays $[x^m, x^{\nu})$.
The intersection points are just given by  $ x(H,\nu ), 1\leq \nu\leq m-1$.

Let $x_j(H,m)$ be the $j$-th coordinate of the point $x(H,m)$.
 Now
$$
\left| x_j(H,\nu )-\frac{L}{m} \right|=
\begin{cases}
\frac{H}{m\nu}, \,\,\, 1\le j\le\nu, 
\\
\frac{H}{m(m-\nu)}, \,\,\, \nu+1 \le j\le m, 
\end{cases}
$$
and
$$
\max_{1\le j,\nu \le m} \left| x_j(H,\nu )-\frac{L}{m} \right| \le H.
$$

This  means that $\Omega_H\subset \Omega_H^*\subset \left[\frac{L}{m} -H, \frac{L}{m}\right]^m$ and so (\ref{b}) follows.
\end{proof}

Since $\rho_H \le M_H$, 
we see that
$$
N_n^{(m)}(L) \le c(q,m) \sum_{H=0}^{(m-1)L}\frac{(H+1)^m}{q^H} q^{2mn -(m+1)L} \le c_1(q,m)
 q^{2mn -(m+1)L},
$$
where
$$
c_1(q,m):= 
c(q,m) \times \sum_{H=0}^{\infty} \frac{(H+1)^m}{q^H} .
$$
 The last estimate together with Lemma 2 leads to
\begin{equation*}
N_n^{(m)}(L) \le c_2(q,m)
 q^{\min \{(m+1)L,2mn -(m+1)L\}}
\end{equation*}
with $c_2(q,m)= \max ( c_1(q,m), 1)$.

This is very close to the upper bound
of Theorem 24 in \cite{MM}.

As an immediate consequence, we have  the inequality
\begin{equation*}
N_n^{(m)}(L)
\ll
q^{nm}\times
q^{-|(m+1) L - mn | }.
\end{equation*}

\newpage 
{\bf Proof of Theorem 1.}

For fixed $n$ and $m$, we 
introduce the linear complexity deviation
$\Delta_L:= \Delta_n^{(m)}(L) = L-\left\lceil\frac{nm}{m+1}\right\rceil $. Then
the number $Z_n^{(m)}(\Delta)$ of sequences with linear complexity
deviation $\Delta$ satisfies
\begin{equation*}
Z_n^{(m)}(\Delta)
\ll
q^{nm}\times
q^{-|\Delta|(m+1)}.
\end{equation*}

Now
\[E_n^{(m)} =  q^{-nm}
\sum_{L=0}^n
L\cdot
N_n^{(m)}(L)
=
\left\lceil\frac{mn}{m+1}\right\rceil
+
q^{-nm}
\sum_{L=0}^n
\Delta_L\cdot
N_n^{(m)}(L)
\]
\[
=
\left\lceil\frac{mn}{m+1}\right\rceil
+
q^{-nm}
\sum_{\Delta=-\infty}^\infty
\Delta\cdot
Z_n^{(m)}(\Delta)
=
\left\lceil\frac{mn}{m+1}\right\rceil+O(1)
\]

The proof is complete.\hfill $\Box$

\begin{conjecture}
In view of the numerical results in 
{\rm \cite{MM}}
we conjecture
$$
\left\lceil\frac{nm}{m+1}\right\rceil-1+O(1/n)
\le
E_n^{(m)}
\le\left\lceil\frac{nm}{m+1}\right\rceil+1+O(1/n).$$
\end{conjecture}

\end{document}